\documentclass[11pt]{article}
\usepackage{IEEEtrantools}
\usepackage{mathtools, nccmath}
\usepackage{latexsym}
\usepackage{amsfonts}
\usepackage{amssymb}
\usepackage{amsmath,amsfonts,amsthm}
\usepackage{mathrsfs}
\usepackage{enumerate}
\usepackage{dsfont}
\usepackage{hyperref}
\usepackage{graphicx}
\usepackage{epstopdf}
\usepackage{tikz}
\usetikzlibrary{shapes,arrows.meta,positioning}
\usepackage{cleveref}
\usepackage{tikz-cd}
\usepackage{psfrag}

\usepackage{multirow}
\usepackage{tikz}
\usetikzlibrary{arrows,decorations.markings}
\usepackage{bookmark}
\usepackage{hyperref}
\hypersetup{
     colorlinks   = true,
     citecolor    = blue
}

\newcommand{\newsection}[1]{\setcounter{equation}{0}
\setcounter{dfn}{0}
\section{#1}}

\newtheorem{dfn}{Definition}[section]
\newtheorem{thm}[dfn]{Theorem}
\newtheorem{lmma}[dfn]{Lemma}
\newtheorem{ppsn}[dfn]{Proposition}

\newtheorem{rmrk}[dfn]{Remark}

\newcommand{\bbc}{\mathbb{C}}
\newcommand{\bbn}{\mathbb{N}}

\def \qed { \mbox{}\hfill
$\Box$\vspace{1ex}}

\title{Sections and Chapters}

\begin{document}

\author{\sc{Keshab Chandra Bakshi\,,\,Satyajit Guin\,,\,Debabrata Jana}\,}
\title{A few remarks on intermediate subalgebras of an inclusion of $C^*$-algebras}
\maketitle


\begin{abstract}
We show that the angle between intermediate $C^*$-subalgebras of an inclusion of simple $C^*$-algebras with finite Watatani index is stable. The notion of angle is instrumental in providing a bound for the cardinality of the lattice of intermediate subalgebras for an irreducible inclusion of simple $C^*$-algebras with finite Watatani index. We improve the existing upper bound for the cardinality of this set.
\end{abstract}
\bigskip

{\bf AMS Subject Classification No.:} {\large 46}L{\large 05}, {\large 47}L{\large 40}.
\smallskip


\hypersetup{linkcolor=blue}
\bigskip

\bigskip

\newsection{Introduction}\label{Sec 0}

The notion of the index $[M : N]$ for a subfactor $N\subset M$ of type $II_1$ factors was introduced by Jones in his seminal paper \cite{J}, along with the notion of the basic construction. Later Kosaki \cite{Ko} generalized the notion of index and basic construction in terms of suitable conditional expectations for subfactors of any type. Inclusions of simple $C^*$-algebras encompass both the type $II_1$ and ($\sigma$-finite) type $III$ subfactor theory. Motivated by the Pimsner-Popa basis \cite{PP}, Watatani \cite{W} generalized Jones’ and Kosaki’s indices to the index of a conditional expectation associated with a unital inclusion of $C^*$-algebras. In the same article, using the language of Hilbert $C^*$-module, Watatani also proposed an analogous notion of basic construction for any pair of unital inclusions of $C^*$-algebras with respect to a finite-index conditional expectation. Watatani's $C^*$-index theory has become an active area of research. Due to the recent breakthrough results in the $C^*$-algebra classification program, it seems to be the right time to explore considerably the inclusion of $C^*$-algebras, their symmetries, and $C^*$-index theory, with the subfactor theory as the principal guide (see \cite{BG,BGS,R}, for instance).

In subfactor theory, the lattice consisting of all intermediate subfactors has attracted much attention over the years, and a substantial amount of work has been done in this direction. Indeed, Watatani had proved that given an irreducible subfactor of type $II_1$ factor with a finite Jones index, the lattice of intermediate subfactors forms a finite set. Subsequently, Teruya and Watatani \cite{TW} proved that the same result also holds for type $III$-subfactors (with Jones' index replaced by Kosaki's index). Finally, to put both the type $II_1$ and type $III$ factors under the single umbrella, Ino and Watatani proved that the parallel event holds in the $C^*$-world also. Indeed, if we consider an irreducible inclusion of simple $C^*$-algebras with `finite Watatani index', the lattice of all intermediate subalgebras is also a finite set. Given a subfactor $N\subset M$ (both type $II_1$ and type $III$) with $N^{\prime}\cap M=\mathbb{C}$ and $[M:N]<\infty$, Longo \cite{Lon} proved that the cardinality of the lattice of intermediate subfactors $\mathcal{L}(N\subset M)$ is bounded above by ${\big({[M:N]}^2\big)}^{{[M:N]}^2}$, and asked whether it can be bounded by ${[M:N]}^{[M:N]}$. In \cite{BDLR}, this question is answered for the case of type $II_1$ factors by proving that $\texttt{\#}\,\mathcal{L}({N\subset M})\leq 9^{[M:N]}$, and for the case of type $III$ factors, the same result holds \cite{BG}. Furthermore, in \cite{BG} it is shown that the cardinality of the lattice $\mathcal{I}(\mathcal{B\subset A})$ of intermediate subalgebras of an irreducible inclusion of simple $C^*$-algebras $\mathcal{B\subset A}$ with a conditional expectation of index-finite type is bounded above by $9^{[\mathcal{A:B}]^{\,2}_0}$, where $[\mathcal{A:B}]_0$ is the so-called `minimal index'.

To obtain the upper bound for the case of type $II_1$ factors, in \cite{BDLR} a new notion of `angle' between intermediate subfactors is introduced, and a crucial rigidity phenomenon of the angle is observed that shows a surprising connection between the angle and the kissing number in geometry. On the other hand, the main novelty of \cite{BG} is the introduction of the analogous notion of angle between intermediate $C^*$-subalgebras of $\mathcal{B\subset A}$. A similar rigidity result also holds for angle in the $C^*$ set-up. In fact, this rigidity of the angle is exploited to obtain the upper bound of the cardinality of the set $\mathcal{I}(\mathcal{B\subset A})$.

In this paper, our goal is twofold. On the one hand, we improve the bound for the cardinality of the set $\mathcal{I}(\mathcal{B\subset A})$ by proving that $\texttt{\#}\mathcal{I}(\mathcal{B\subset A})\leq 9^{[\mathcal{A:B}]_0}$. We remark that as both the type $II_1$ factor and ($\sigma$-finite) type $III$ factors are particular examples of simple $C^*$-algebras, this result can be thought of as the apt generalization of the corresponding results in \cite{BDLR,BG} mentioned earlier. On the other hand, we prove the stability of the angle between intermediate $C^*$-subalgebras, that is, angle remains invariant under tensor product with the $C^*$-algebra of compact operators on a separable Hilbert space.


\newsection{Preliminaries}

\subsection{Finite-index conditional expectations}

Motivated by Jones' index theory \cite{J} for subfactors and Pimsner-Popa basis \cite{PP}, Watatani developed a theory of index for inclusion of $C^*$-algebras that generalizes the Jones index for a subfactor of type $II_1$ factor, and also Kosaki's index \cite{Ko} for a subfactor of type $III$ factor.  Given a pair $\mathcal{B \subset A}$ of $C^*$-algebras, a conditional expectation $E:\mathcal{A\to B}$ is said to be of \textit{index-finite type} if there exists a finite set $\{\lambda_1,\ldots,\lambda_n\}\subset\mathcal{A}$ such that
\begin{center}
$x=\sum_{i=1}^nE(x\lambda_i)\lambda^*_i=\sum_{i=1}^n\lambda_iE(\lambda^*_ix)$
\end{center}
for every $x\in\mathcal{A}$ \cite{Wa}. Such a set $\{\lambda_1,\ldots,\lambda_n\}$ is called a \textit{quasi-basis} for $E$ and the Watatani index of $E$ is defined by
\begin{center}
$\mathrm{Ind}_w(E):= \sum_{i=1}^n \lambda_i\lambda^*_i$.
\end{center}
It is known that $\mathrm{Ind}_w(E)$ is not a scalar but a positive invertible element in $\mathcal{Z(A)}$, and is independent of the quasi-basis $\{\lambda_i\}$. In particular, if $\mathcal{A}$ is a simple $C^*$-algebra, the index is scalar-valued. 

We denote by $\mathcal{E}_0(\mathcal{A,B})$ the set of all index-finite type conditional expectations from $\mathcal{A}$ onto $\mathcal{B}$. A conditional expectation $F\in\mathcal{E}_0(\mathcal{A , B})$ is said to be minimal if it satisfies $\mathrm{Ind}_w(F) \leq \mathrm{Ind}_w(E)$ for all $E\in\mathcal{E}_0(\mathcal{A, B})$ (see \cite{Wa}). For inclusion of simple $C^*$-algebras, we have a privileged minimal conditional expectation as mentioned below.

\begin{thm}{\em \cite[Theorem $2.12.3$]{Wa}} Let $\mathcal{B\subset A}$ be an inclusion of simple $C^*$-algebras such that $\,\mathcal{E}_0(\mathcal{A,B})\neq\emptyset$. Then, there exists a unique minimal conditional expectation $E_0$ from $\mathcal{A}$ onto $\mathcal{B}$.
\end{thm}

For inclusion of simple $C^*$-algebras $\mathcal{B\subset A}$, the minimal index is defined as
\begin{center}
$[\mathcal{A:B}]_0:=\mathrm{Ind}_w(E_0).$
\end{center}
This seems to be a proper place to point out that if $N\subset M$ is a subfactor with finite Jones index $[M:N]$ and is irreducible (i.e., $N^{\prime}\cap M=\mathbb{C}$), then the trace preserving conditional expectation $E^M_N$ is the minimal conditional expectation with $[M:N]={[M:N]}_0$. In general, for reducible subfactors, the minimal index and Jones index need not coincide.

\subsection{Watatani’s $C^*$-basic construction}

In the subfactor theory \cite{JS}, Jones' basic construction plays a pivotal role. Using the language of the Hilbert $C^*$-module, Watatani proposed a parallel notion of basic construction in the $C^*$-world, the so-called $C^*$-basic construction \cite{Wa}. Let $\mathcal{B\subset A}$ be an inclusion of unital $C^*$-algebras and $E:\mathcal{A\to B}$ be a faithful conditional expectation. Then, $\mathcal{A}$ is a pre-Hilbert right $\mathcal{B}$-module with respect to the $\mathcal{B}$-valued inner product given by
\begin{equation}\label{B-valued}
\langle x, y\rangle_\mathcal{B}=E_\mathcal{B}(x^*y)\quad \text{ for all } x, y\in\mathcal{A},
\end{equation}
and $\mathscr{U}$ denotes the Hilbert (right) $\mathcal{B}$-module completion of $\mathcal{A}$. The space of adjointable maps on $\mathscr{U}$, denoted by $\mathcal{L}_\mathcal{B}(\mathscr{U})$, is a unital $C^*$-algebra and $\mathcal{A}$ embeds in it as a unital $C^*$-subalgebra. There exists a projection $e_\mathcal{B}\in\mathcal{L}_\mathcal{B}(\mathscr{U})$ (called the Jones projection associated to $E$) such that $e_\mathcal{B}ae_\mathcal{B}=E(a)e_\mathcal{B}$ for all $a\in\mathcal{A}$. Consider $\mathcal{A}_1:=\overline{\mbox{span}}\{xe_\mathcal{B}y:x,y\in\mathcal{A}\}\subseteq\mathcal{L}_\mathcal{B}(\mathscr{U})$, which turns out to be a $C^*$-algebra (not always unital) and is called the $C^*$-basic construction of the inclusion $\mathcal{B\subset A}$.

\begin{rmrk}\rm
If $E:\mathcal{A\to B}$ has finite index with a quasi-basis $\{\lambda_i\}$, then we have the following\,:
\begin{enumerate}[$(i)$]
\item The two norms $||.||_\mathcal{A}$ and $||.||$ on $\mathcal{A}$, where $||.||$ is the given $C^*$-norm and $||x||_\mathcal{A}:=||E_\mathcal{B}(x^*x)||^{\frac{1}{2}}$, are equivalent (see \cite{Wa}, or Lemma 2.11 in \cite{BG}). In particular, $\mathcal{A}$ itself is a Hilbert $\mathcal{B}$-module.
\item $\mathcal{A}_1$ is unital and equal to $C^*(\mathcal{A},e_\mathcal{B})$ (see Proposition 2.1.5 in \cite{Wa}).
\item There exists a finite-index conditional expectation $E_1:\mathcal{A}_1\to\mathcal{A}$ (called the dual conditional expectation) with a quasi-basis $\{\lambda_ie_B(\mbox{Ind}_w(E))^\frac{1}{2}\}$ which satisfies $E_1(xe_\mathcal{B}y)=\mbox{Ind}_w(E)^{-1}xy$ for all $x,y\in\mathcal{A}$ and $\mbox{Ind}_w(E_1)=\sum_i\lambda_iE(\mbox{Ind}_w(E))e_\mathcal{B}\lambda_i^*$. Moreover, if $\mbox{Ind}_w(E)\in\mathcal{B}$, then $\mbox{Ind}_w(E_1)=\mbox{Ind}_w(E)$ (see Proposition 2.3.2 and 2.3.4 in \cite{Wa}).
\item $\mathcal{A}_1=\mbox{span}\{xe_\mathcal{B}y:x,y\in\mathcal{A}\}=:\mathcal{A}e_\mathcal{B}\mathcal{A}$ (see \cite{Wa}, Lemma 2.2.2).
\end{enumerate}
\end{rmrk}

\noindent\textbf{Non-unital situation\,:} For the case of inclusions $\mathcal{B\subset A}$ of non-unital $C^*$-algebras, there is no way to define index in terms of quasi-basis, although the basic constuction makes sense. This is because a quasi-basis exists only if both $\mathcal{A}\mbox{ and }\mathcal{B}$ are unital \cite{Wa}. It is also easy to show that a simple unitization trick does not resolve this problem, and we invite the reader to visit Section $2.3$ in \cite{Iz}.

 Izumi's aim in \cite{Iz} was to supplement Watatani's $C^*$-basic construction by using a second dual argument in the non-unital case, and define the index `$\mbox{Ind}_w(E)$' for some class of inclusions (namely, for which $\mathcal{A}\subset\mathcal{A}_1$) of non-unital $C^*$-algebras. First consider the multiplier algebras $M(\mathcal{A})\mbox{ and }M(\mathcal{B})$ as $C^*$-subalgebras of $\mathcal{A}^{**}\mbox{ and }\mathcal{B}^{**}$ respectively. When $\mathcal{AB}=\mathcal{A}$ holds, we have $M(\mathcal{B})\subset M(\mathcal{A})$ (in general, no such inclusion holds). Recall the following facts from Lemma $2.6$, \cite{Iz}\,:
\begin{enumerate}[$(i)$]
\item If $\mbox{Ind}_pE<\infty$, we have $\mathcal{A}=\mathcal{AB}$. In consequence, $M(\mathcal{B})$ is a unital subalgebra of $M(\mathcal{A})$;
\item The restriction of $E^{**}$ to $M(\mathcal{A})$ is a conditional expectation from $M(\mathcal{A})$ onto $M(\mathcal{B})$.
\end{enumerate}
Next, with the help of the bounded normal operator valued weight $\widehat{E^{**}}$, Izumi defines the index of $E^{**}$ as $\widehat{E^{**}}(1)\in\mathcal{Z}(\mathcal{A}^{**})$. Finally, one has the following.
\begin{dfn}[Def. $2.10$, \cite{Iz}]
For $\mathcal{B}\subset^{\,E}\mathcal{A}$, if $\mathcal{A}\subset\mathcal{A}_1$ holds, the Watatani index of $E$ is defined by $\mbox{Ind}_w(E):=\widehat{E^{**}}(1)\in\mathcal{Z}(M(\mathcal{A}))$.  The dual conditional expectation $E_1:\mathcal{A}_1\to\mathcal{A}$ is defined by $E_1(x):=(\mbox{Ind}_wE)^{-1}\widehat{E}(x)$, and $E_1$ satisfies $\text{Ind}_pE_1\leq||\mbox{Ind}_wE||$. When $\mathcal{A}\subset\mathcal{A}_1$ does not hold, $\mbox{Ind}_wE:=\infty$.
\end{dfn}
\noindent Note that $\mathcal{A}\subset\mathcal{A}_1$ holds if and only if $\widehat{E^{**}}(1)\in\mathcal{Z}(M(\mathcal{A}))$ (Thm. $2.8$, \cite{Iz}).
\smallskip

For example, given a unital inclusion $\mathcal{B\subset A}$, consider the corresponding non-unital inclusion $\mathcal{B}\otimes K(\mathcal{H})\subset \mathcal{A}\otimes K(\mathcal{H})$, where $\mathcal{H}$ is a separable Hilbert space and $K(\mathcal{H})$ denotes the space of compact operators. In this case, if $E:\mathcal{A\to B}$ is of finite index, then $E\otimes\mbox{id}$ is a conditional expectation of index-finite type with $\mbox{Ind}_w(E\otimes\mbox{id})=\mbox{Ind}_w(E)\otimes 1$ (see Lemma $2.11$ in \cite{Iz}).

\subsection{Intermediate $C^*$-subalgebras and angles}

Let $\mathcal{B\subset A}$ be an inclusion of unital $C^*$-algebras with a conditional expectation $E:\mathcal{A\to B}$. Then, for the $\mathcal{B}$-valued inner product $\langle .,.\rangle_\mathcal{B}$ on $\mathcal{A}$ given by $\langle x,y\rangle_\mathcal{B}=E(x^*y)$, one has the following well known analogue of the Cauchy-Schwarz inequality
\begin{eqnarray}\label{an inequality}
||\langle x,y\rangle_\mathcal{B}||\leq ||x||_\mathcal{A}||y||_\mathcal{A}
\end{eqnarray}
for all $x,y\in\mathcal{A}$, where $||x||_\mathcal{A}:=||E_\mathcal{B}(x^*x)||^{1/2}$. If $\mathcal{B\subset A}$ is irreducible, that is $\mathcal{B}^{\prime}\cap \mathcal{A}=\mathbb{C}$, then any intermediate $C^*$-subalgebra $\mathcal{C}$ is simple and the conditional expectation $F:\mathcal{A\to C}$ satisfies the compatibility condition $E=E|_\mathcal{C}\circ F$ (see \cite{Iz}). Employing \Cref{an inequality}, the notion of interior and exterior angles has been introduced in \cite{BG}, which further generalizes the angle between a pair of intermediate subfactors of a given finite index subfactor, first introduced in \cite{BDLR}. A surprising connection between angle and the kissing number \cite{PZ} in geometry/sphere packing has been discovered.

With the help of \cite{Iz}, one can easily generalize angle in the non-unital case as follows. Recall that if $\mathcal{B\subset A}$ is an inclusion of $\sigma$-unital infinite dimensional $C^*$-algebras with a conditional expectation $E$ from $\mathcal{A}$ onto $\mathcal{B}$ with $\mathcal{A}$ simple, then $\text{Ind}_p(E)=\text{Ind}_w(E)$ (Corollary $3.7$, \cite{Iz}), where $\text{Ind}_p E$ is the Pimsner-Popa probabilistic index (\cite{PP}, or Definition $2.1$ in \cite{Iz}) defined as follows
\[
\text{Ind}_p E:=\text{inf}\big\{\lambda>0\,:\,\mathrm{Id}-\frac{1}{\lambda} E~\text{ is positive}\big\}.
\]
\begin{dfn}[Def. $5.1$, \cite{BG}]\label{angle}
Consider an inclusion of $\sigma$-unital simple $C^*$-algebras $\mathcal{B \subset A}$ with a conditional expectation $E$ satisfying $\mathrm{Ind}_w(E)<\infty$. If $\mathcal{B\subset A}$ is irreducible, then for a pair of intermediate $C^*$-subalgebras $\mathcal{C}$ and $\mathcal{D}$, we define the angle $\alpha^{\mathcal{B}}_{\mathcal{A}}(\mathcal{C,D})$ by the following expression
\[ 
\cos\left(\alpha^{\mathcal{B}}_{\mathcal{A}}(\mathcal{C,D})\right)=\displaystyle \frac{\lVert{\langle
    e_{\mathcal{C}}-e_{\mathcal{B}},e_{\mathcal{D}}-e_{\mathcal{B}}\rangle}_{\mathcal{A}}\rVert}{{\lVert
    e_\mathcal{C}-e_\mathcal{B}\rVert}_{\mathcal{A}}{\lVert
    e_\mathcal{D}-e_\mathcal{B}\rVert}_{\mathcal{A}}}.
\]
\end{dfn}

By definition, notice that the angle is allowed to take values only in the interval $[0,\pi/2]$. 


\newsection{Main results}

We prove two results in this article. Firstly, for an irreducible unital inclusion $\mathscr{B\subset A}$ of simple $C^*$-algebras with a conditional expectation of index-finite type, the angle in \Cref{angle} remains invariant under tensor product with the $C^*$-algebra of compact operators on a separable Hilbert space. Secondly, for such an inclusion, the bound for the cardinality of the lattice of intermediate $C^*$-subalgebras is $9^{[\mathscr{A:B}]_0}$, where $[\mathscr{A:B}]_0$ is the minimal index.

\subsection{Angle between intermediate subalgebras and its stability}

Aim of this subsection is to prove the following theorem.

\begin{thm}\label{anglestable}
Let $\mathcal{H}$ be a separable Hilbert space and $K(\mathcal{H})$ be the $C^*$-algebra of compact operators on $\mathcal{H}$. For an irreducible unital inclusion $\mathscr{B\subset A}$ of simple $C^*$-algebras with a conditional expectation of index-finite type and intermediate $C^*$-subalgebras $\mathscr{B\subset C,D\subset A}$, the angle is stable under tensor product with $K(\mathcal{H})$. That is, \[\alpha^{\mathscr{B}}_{\mathscr{A}}(\mathscr{C,D})= \alpha^{\mathscr{B}\otimes K(\mathcal{H})}_{\mathscr{A}\otimes K(\mathcal{H})}\big(\mathscr{C}\otimes K(\mathcal{H}),\mathscr{D}\otimes K(\mathcal{H})\big).\]
\end{thm}

We first need a few preliminary results. Let $\mathcal{I}(\mathscr{B}\subset \mathscr{A})$ be the lattice of intermediate  $C^*$-subalgebras. We first need the following result that is related to the tensor spilitting theorem in the unital case proved in \cite{Z, Zs}.

\begin{lmma}\label{ror}
Let $\mathscr{B\subset A}$ be a unital inclusion of $C^*$-algebras. There is a bijective correspondence between the sets $\mathcal{I}\big(\mathscr{B}\otimes K(\mathcal{H})\subset \mathscr{A}\otimes K(\mathcal{H})\big)$ and $\mathcal{I}(\mathscr{B}\subset \mathscr{A})$.
\end{lmma}
\begin{prf}
Fix an orthonormal basis $\{e_i:i=0,1,\ldots\}$ of $\mathcal{H}$ and take $p=|e_0\rangle\langle e_0|$, the rank one projection onto $\bbc e_0$. Also, denote $p_n=|e_n\rangle\langle e_n|$, the rank one projection onto $\bbc e_n$, for all $n\in\bbn$. Let $\ell$ denote the right shift operator on $\mathcal{H}$. For any intermediate $C^*$-subalgebra $\mathscr{B}\otimes K(\mathcal{H})\subset\mathcal{M}\subset \mathscr{A}\otimes K(\mathcal{H})$, consider $\mathscr{C}=\{a\in \mathscr{A}:a\otimes p\in\mathcal{M}\}$. Clearly $\mathscr{C}$ is a unital $C^*$-subalgebra of $\mathscr{A}$ and we have $\mathscr{B\subset C\subset A}$. We claim that $\mathscr{C}\otimes K(\mathcal{H})=\mathcal{M}$.

Choose $\xi=\sum a_j\otimes T_j\in \mathscr{C}\otimes K(\mathcal{H})$, where $T_j$'s are finite rank operators. For each $j$, since $a_j\in\mathscr{C}$, by definition we have $a_j\otimes p\in \mathcal{M}$. Since $\mathscr{B}\otimes K(\mathcal{H})\subset\mathcal{M}$, we get that $a_j\otimes p_n=(1\otimes \ell^np)(a_j\otimes p)(1\otimes p\ell^{-n})\in\mathcal{M}$ for any $n\in\mathbb{N}$, where $\ell^{-n}=(\ell^*)^n$. Since $T_j$'s are finite rank operators, we get $a_j\otimes T_j\in\mathcal{M}$ for each $j$, and consequently, $\mathscr{C}\otimes K(\mathcal{H})\subset\mathcal{M}$ by density of finite rank operators in $K(\mathcal{H})$ and closedness of $\mathcal{M}$.
 
Now, to show that $\mathcal{M}$ is contained in $\mathscr{C}\otimes K(\mathcal{H})$, let $\{E_{ij}:i,j\in\bbn\cup\{0\}\}$ be the matrix unit for $K(\mathcal{H})$. Note that $1\otimes K(\mathcal{H})$ is contained in $\mathscr{B}\otimes K(\mathcal{H})$, and hence in $\mathcal{M}$. Set $P_n = \sum_{i=0}^n 1\otimes E_{ii}$. Then, $||x – P_n x P_n||\to 0$ for $n\to\infty$ for all $x\in\mathscr{A}\otimes K(\mathcal{H})$, and hence for all $x\in\mathcal{M}$. Thus, given $x\in\mathcal{M}$, it suffices to show that $P_nxP_n$ belongs to $\mathscr{C}\otimes K(\mathcal{H})$ for all $n$. Now, $P_nxP_n=\sum_{i,j=0}^nx_{ij}\otimes E_{ij}$, and $x_{ij}\otimes E_{ij}=(1\otimes E_{ii})x(1\otimes E_{jj})$ belongs to $\mathcal{M}$. As $x_{ij} \otimes p=x_{ij} \otimes E_{00}=(1\otimes E_{0i})(x_{ij}\otimes E_{ij})(1\otimes E_{j0})$ belongs to $\mathcal{M}$, we conclude that $x_{ij}$ belongs to $\mathscr{C}$. This shows that $P_nxP_n$ belongs to $\mathscr{C}\otimes K(\mathcal{H})$, and consequently $x\in\mathscr{C}\otimes K(\mathcal{H})$.

Therefore, any $\mathscr{B}\otimes K(\mathcal{H})\subset \mathcal{M}\subset \mathscr{A}\otimes K(\mathcal{H})$ is of the form $\mathscr{C}\otimes K(\mathcal{H})$ for some $\mathscr{B\subset C\subset A}$, and consequently the map $\phi:\mathcal{M}\mapsto\mathscr{C}$ becomes bijective.\qed
\end{prf}

We sincerely thank Professor Mikael Rørdam for helping us with the proof of the above result and kindly allowing us to present it here. Note that in the above proof neither the simpleness nor the index finite-type are needed for the inclusion $\mathscr{B\subset A}$.
\begin{lmma}\label{need}
For any unital inclusion $\mathscr{B\subset A}$ of $C^*$-algebras, we have $M(\mathscr{B}\otimes K(\mathcal{H}))\subset M(\mathscr{A}\otimes K(\mathcal{H}))$, where $M(.)$ denotes the multiplier algebra.
\end{lmma}
\begin{prf}
Note that for $M(\mathscr{B}\otimes K(\mathcal{H}))\subset M(\mathscr{A}\otimes K(\mathcal{H}))$, it is enough to show that (see Section $2.3$, \cite{Iz}) $\mathscr{A}\otimes K(\mathcal{H})=(\mathscr{A}\otimes K(\mathcal{H}))(\mathscr{B}\otimes K(\mathcal{H}))=\mathscr{AB}\otimes K(\mathcal{H})K(\mathcal{H})=\mathscr{A}\otimes K(\mathcal{H})K(\mathcal{H})$, where the last equality follows because $\mathscr{B\subset A}$ is a unital inclusion and hence we have $\mathscr{A=AB}$. Since, $K(\mathcal{H})K(\mathcal{H})\subseteq K(\mathcal{H})$, it is further enough to show that $K(\mathcal{H})\subseteq K(\mathcal{H})K(\mathcal{H})$. However, this is easy as for any $T\in K(\mathcal{H})$, first write $T=U|T|$ by the polar decomposition. Now, $|T|$ is positive and compact, and hence so is its positive square root $S$. Then, $T=(US)S$ shows that $K(\mathcal{H})\subseteq K(\mathcal{H})K(\mathcal{H})$.\qed
\end{prf}

Due to \Cref{need}, the set-up in Section $2.3$, \cite{Iz}, for the Watatani index in the non-unital situation is applicable in our situation. If $\mathscr{B\subset A}$ is a unital inclusion of simple $C^*$-algebras with finite Watatani index, then the inclusion $\mathscr{B}\otimes K(\mathcal{H})\subset\mathscr{A}\otimes K(\mathcal{H})$ is also so. This follows from Lemma $2.11$ in \cite{Iz} (observe that unitality and simpleness of $\mathscr{A}$ gives us $\mathcal{Z}(M(\mathscr{A}))=\mathcal{Z}(\mathscr{A})=\bbc$).
\smallskip

The following result is expected to be known, but unfortunately, we could not find a proof in the literature.

\begin{ppsn}\label{basic}
Let $\mathscr{B\subset A}$ be a finite-index unital $C^*$-inclusion. The basic construction of $\mathscr{B}\otimes K(\mathcal{H})\subset \mathscr{A}\otimes K(\mathcal{H})$ is $\mathscr{A}_1\otimes K(\mathcal{H})$, where $\mathscr{B\subset A\subset A}_1$ is the $C^*$-basic construction.
\end{ppsn}
\begin{prf}
Let $\mathscr{B\subset A}$ be a finite-index unital $C^*$-inclusion with respect to conditional expectation $E$ and $\mathscr{B\subset A\subset A}_1$ be the $C^*$-basic construction with Jones projection $e_E\in\mathscr{A}_1$. Suppose that $\widetilde{\mathscr{C}}$ is the $C^*$-basic construction of $\mathscr{B}\otimes K(\mathcal{H})\subset \mathscr{A}\otimes K(\mathcal{H})$. By Section $2.3$ in \cite{Iz} (see just before Thm. $2.8$), we know that if $E:\mathscr{A}\to\mathscr{B}$ is the conditional expectation with Jones projection $e_E$, then the Jones projection $e_{E^{**}}$ corresponding to $E^{**}:\mathscr{A}^{**}\to\mathscr{B}^{**}$ can be identified with $e_E$. Since $\big(\mathscr{B}\otimes K(\mathcal{H})\big)^{**}=\mathscr{B}^{**}\otimes B(\mathcal{H})$ (and similarly, $\mathscr{A}$ in place of $\mathscr{B}$) and $E^{**}\otimes\mbox{id}:\mathscr{A}^{**}\otimes B(\mathcal{H})\to\mathscr{B}^{**}\otimes B(\mathcal{H})$ is the conditional expectation with Jones projection $e_{E^{**}}\otimes\mbox{id}$, we see that the Jones projection for the inclusion $\mathscr{B}\otimes K(\mathcal{H})\subset\mathscr{A}\otimes K(\mathcal{H})$ can be identified with $e\otimes\mbox{id}$. Since $\widetilde{\mathscr{C}}$ is the $C^*$-subalgebra of the unital $C^*$-algebra $\mathcal{L}_{\mathscr{B}\otimes K(\mathcal{H})}\big(\mathscr{A}\otimes K(\mathcal{H})\big)$ defined by the following
\[
\widetilde{\mathscr{C}}:=\overline{\mathrm{span}}\{(a_1\otimes T_1)(e\otimes\mbox{id})(a_2\otimes T_2):a_1,a_2\in\mathscr{A};\,T_1,T_2\in K(\mathcal{H})\},
\]
we see that inside the ambient unital $C^*$-algebra $\mathscr{A}_1\otimes B(\mathcal{H}),\,\widetilde{\mathscr{C}}$ can be viewed as a subalgebra of $\mathscr{A}_1\otimes K(\mathcal{H})$, where $\mathscr{A}_1=\mathscr{A}e\mathscr{A}$. Since $\mathscr{A}\otimes K(\mathcal{H})\subset\widetilde{\mathscr{C}}\subseteq\mathscr{A}_1\otimes K(\mathcal{H})$, by \Cref{ror} we have $\widetilde{\mathscr{C}}=\mathscr{C}\otimes K(\mathcal{H})$ for some unital $C^*$-subalgebra $\mathscr{A\subset C\subseteq A}_1$. We claim that $\mathscr{C}=\mathscr{A}_1$. For any minimal projection $p\in K(\mathcal{H})$ (in particular, $p=|e_0\rangle\langle e_0|$), we have $(1\otimes p)(e\otimes\mbox{id})(1\otimes p)=e\otimes p\in\widetilde{\mathscr{C}}=\mathscr{C}\otimes K(\mathcal{H})$. Choose a functional $\psi_p\in K(\mathcal{H})^*$ such that $\psi_p(p)=1$, and consider $\mbox{id}\otimes\psi_p:\mathscr{C}\otimes K(\mathcal{H})\to \mathscr{C}$. Then, $(\mbox{id}\otimes\psi_p)(e\otimes p)=e\in\mathscr{C}$. Since $\mathscr{A}_1$ is generated by $\mathscr{A}$ and $e$, and we already have $\mathscr{A}\subset\mathscr{C}$, we get that $\mathscr{A}_1\subseteq\mathscr{C}$, and consequently $\widetilde{\mathscr{C}}=\mathscr{C}\otimes K(\mathcal{H})=\mathscr{A}_1\otimes K(\mathcal{H})$, which completes the proof.\qed
\end{prf}

By the previous \Cref{basic}, the tower of basic construction for the inclusion $\mathscr{B}\otimes K(\mathcal{H})\subset \mathscr{A}\otimes K(\mathcal{H})$ is the following
\[
\mathscr{B}\otimes K(\mathcal{H})\subset \mathscr{A}\otimes K(\mathcal{H})\subset \mathscr{A}_1\otimes K(\mathcal{H})\subset\mathscr{A}_2\otimes K(\mathcal{H})\subset\cdots
\]
where $\mathscr{B\subset A\subset A_1\subset A}_2\subset\cdots$ is the tower of basic construction for the unital inclusion $\mathscr{B}\subset\mathscr{A}$. Following \cite{Iz}, for $k\geq 0$, we define $(\mathscr{B}\otimes K(\mathcal{H}))^\prime\cap (\mathscr{A}_k\otimes K(\mathcal{H})):=(\mathscr{B}\otimes K(\mathcal{H}))^\prime\cap M(\mathscr{A}_k\otimes K(\mathcal{H}))$. This definition is consistent due to \Cref{need}. Then, we have the following
\begin{IEEEeqnarray}{lCl}\label{prime}
(\mathscr{B}\otimes K(\mathcal{H}))^\prime\cap (\mathscr{A}_k\otimes K(\mathcal{H})) &=& (\mathscr{B}\otimes K(\mathcal{H}))^\prime\cap M(\mathscr{A}_k\otimes K(\mathcal{H}))\nonumber\\
&=& (\mathscr{B}^\prime\cap M(\mathscr{A}_k))\otimes\bbc\nonumber\\
&=& \mathscr{B}^\prime\cap \mathscr{A}_k\,,
\end{IEEEeqnarray}
which is a finite-dimensional algebra since $[\mathscr{A:B}]_0<\infty$. Thus, the Fourier theory in \cite{BG} (see also \cite{BGS}) can be extended to this non-unital set-up. Note that by \Cref{prime}, if $\mathscr{B\subset A}$ is irreducible, then so is $\mathscr{B}\otimes K(\mathcal{H})\subset\mathscr{A}\otimes K(\mathcal{H})$.

Now, we are ready to prove \Cref{anglestable}.
\medskip

\textbf{Proof of \Cref{anglestable}\,:} Let $\mathscr{B\subset A}$ be a unital inclusion of simple $C^*$-algebras with a conditional expectation of index-finite type and $e_\mathscr{B}$ be the associated Jones projection. In \Cref{basic}, it is explained that the Jones projection for the inclusion $\mathscr{B}\otimes K(\mathcal{H})\subset\mathscr{A}\otimes K(\mathcal{H})$ is identified with $e_\mathscr{B}\otimes\mbox{id}$. For intermediate $C^*$-subalgebras $\mathscr{C,D}$ of $\mathscr{B\subset A}$, consider $x=e_{\mathscr{C}\otimes K(\mathcal{H})}-e_{\mathscr{B}\otimes K(\mathcal{H})}=(e_\mathscr{C}-e_\mathscr{B})\otimes\mbox{id}$ and $y=e_{\mathscr{D}\otimes K(\mathcal{H})}-e_{\mathscr{B}\otimes K(\mathcal{H})}=(e_\mathscr{D}-e_\mathscr{B})\otimes\mbox{id}$. By \Cref{ror}, we know all intermediate $C^*$-subalgebras of $\mathscr{B}\otimes K(\mathcal{H})\subset\mathscr{A}\otimes K(\mathcal{H})$. Therefore, we have $||\langle x,y\rangle_{\mathscr{A}\otimes K(\mathcal{H})}||=||E_1\big((e_\mathscr{D}-e_\mathscr{B})^*(e_\mathscr{C}-e_\mathscr{B})\big)\otimes\mbox{id}||=||E_1\big((e_\mathscr{D}-e_\mathscr{B})^*(e_\mathscr{C}-e_\mathscr{B})\big)||=||\langle e_\mathscr{C}-e_\mathscr{B},e_\mathscr{D}-e_\mathscr{B}\rangle_\mathscr{A}||$, and similarly we have $||x||_{\mathscr{A}\otimes K(\mathcal{H})}=||e_\mathscr{C}-e_\mathscr{B}||_\mathscr{A}$ and $||y||_{\mathscr{A}\otimes K(\mathcal{H})}=||e_\mathscr{D}-e_\mathscr{B}||_\mathscr{A}$. By \Cref{angle}, the result now follows.\qed



\subsection{On the cardinality of the lattice of intermediate subalgebras }\label{Sec 1}

Aim of this subsection is to prove the following theorem.
\begin{thm}\label{bound}
Let $\mathscr{B\subset A}$ be a unital irreducible inclusion of simple $C^*$-algebras with finite index, that is, $[\mathscr{A:B}]_0<\infty$. Then, the cardinality of the lattice of intermediate $C^*$-subalgebras $\mathcal{I}(\mathscr{B\subset A})$ is bounded by $9^{[\mathscr{A:B}]_0}$.
\end{thm}

Since both the type $II_1$ factor and ($\sigma$-finite) type $III$ factors are particular examples of simple $C^*$-algebras, this result is the apt generalization of the corresponding results in \cite{BDLR,BG}. We describe the literature briefly. Watatani \cite{Wa} showed that if $N \subset M$ is a finite-index irreducible subfactor of type $II_1$, then the lattice $\mathcal{L}(N\subset M)$ of intermediate subfactors is a finite lattice. Subsequently, Teruya and Watatani \cite{TW} showed that $\mathcal{L}(N\subset M)$ is  finite also if $N\subset M$ is a finite-index irreducible subfactor of type $III$. If we consider a unital inclusion $\mathscr{B\subset A}$ of unital $C^*$-algebras, the set of intermediates $C^*$-subalgebras $\mathcal{I}(\mathscr{B\subset A})$ also forms a lattice under the following two operations\,:
\[
\mathscr{C\wedge D}:= \mathscr{C\cap D} ~~\text{and}~~ \mathscr{C\vee D}:= C^*\big(\mathscr{C,D}\big).
\]
Ino and Watatani \cite{Ino} proved that $\mathcal{I}(\mathscr{B\subset A})$ is finite if $\mathscr{A}, \mathscr{B}$ are simple unital $C^*$-algebras with ${\mathscr{B}}^{\prime}\cap\mathscr{A}=\mathbb{C}$ (that is, irreducible) and $[\mathscr{A:B}]_0<\infty$, where $[\mathscr{A:B}]_0$ denotes the minimal index.

On the other hand, Longo \cite{Lon} proved that the number of intermediate subfactors of an irreducible subfactor $N \subset M$ (of any type) with finite-index is bounded by $([M:N]^2)^{[M:N]^2}$ and asked whether this bound could be improved to $[M:N]^{[M:N]}$. Using the notion of angle between intermediate subfactors, authors in \cite{BDLR} proved that for an irreducible subfactor $N \subset M$ of type $II_1$, the bound can be improved significantly to $9^{[M:N]}$. Subsequently, the question for irreducible subfactors of type $III$ factors has been answered in \cite{BG} using the notion of the angle in the $C^*$-settings by proving that the cardinality of the lattice of intermediate subfactors are bounded by $9^{[M:N]}$. Also in \cite{BG}, a bound for the cardinality of the set $\mathcal{I}(\mathscr{B\subset A})$ of all intermediate $C^*$-subalgebras for an inclusion $\mathscr{B\subset A}$ of simple $C^*$-algebras with finite Watatani index is also obtained as $9^{\mathscr{[A:B]}_0^{\,2}}$. In this paper, we further improve this bound as a complete generalization of correponding results for both the type $II_1$ and type $III$ factors.

The following important result is mentioned in \cite{BG} without proof. For the sake of completeness, we start by providing a proof here which is due to Professor Yasuo Watatani through private communication. We sincerely thank him for kindly allowing us to present the proof here.

\begin{ppsn}\label{watatani}
Let $\mathscr{B\subset A}$ be an irreducible $C^*$-inclusion and $\mathscr{C,D}$ be intermediate $C^*$-subalgebras with conditional expectations $E_\mathscr{C}:\mathscr{A}\to\mathscr{C}$ and $E_\mathscr{D}:\mathscr{A}\to\mathscr{D}$ satisfying the compatibility condition. Let $e_\mathscr{C}$ and $e_\mathscr{D}$ are their Jones projections. Then, $e_\mathscr{C}\wedge e_\mathscr{D}=e_{\mathscr{C\cap D}}$.
\end{ppsn}
\begin{prf}
On contrary assume that $e_\mathscr{C}\wedge e_\mathscr{D}\neq e_{\mathscr{C\cap D}}$. Since $e_{\mathscr{C\cap D}}\leq e_\mathscr{C},e_\mathscr{D}$, we have $e_{\mathscr{C\cap D}}\leq e_\mathscr{C}\wedge e_\mathscr{D}$. Let $\mathcal{E}:=\overline{\mathscr{A}}^{\langle\,,\,\rangle_{\mathscr{B}}}$ and $\eta:\mathcal{E}\to\mathscr{A}$ be the canonical inclusion map. There exists some $x\in\mathscr{A}$ such that
\[
\eta(y):=(e_\mathscr{C}\wedge e_\mathscr{D}-e_{\mathscr{C\cap D}})\eta(x)\neq 0\,.
\]
Then, $e_{\mathscr{C\cap D}}\eta(y)=0$ and $(e_\mathscr{C}\wedge e_\mathscr{D})\eta(y)=\eta(y)$. Since $\mathscr{B}^\prime\cap\mathscr{A}_1$ is finite-dimensional and $e_\mathscr{C},e_\mathscr{D}\in\mathscr{B}^\prime\cap\mathscr{A}_1$, we have $(e_\mathscr{C}e_\mathscr{D}e_\mathscr{C})^n$ and $(e_\mathscr{D}e_\mathscr{C}e_\mathscr{D})^n$ converge to $e_\mathscr{C}\wedge e_\mathscr{D}$ in the operator norm topology. Therefore, $(e_\mathscr{C}e_\mathscr{D}e_\mathscr{C})^n\eta(y)$ converges to $(e_\mathscr{C}\wedge e_\mathscr{D})\eta(y)=\eta(y)$ in the Hilbert module norm. Since $(e_\mathscr{C}e_\mathscr{D}e_\mathscr{C})^n\eta(y)=\eta\big((E_\mathscr{C}E_\mathscr{D}E_\mathscr{C})^n(x)\big)\in\eta(\mathscr{C})$, and the Hilbert module norm topology and operator norm topology are identical, we have that $y\in\mathscr{C}$. Similarly, we have $y\in\mathscr{D}$. Therefore, $e_{\mathscr{C\cap D}}\eta(y)=\eta(y)$. This contradicts that $e_{\mathscr{C\cap D}}\eta(y)=0$ and $\eta(y)\neq 0$, which completes the proof.\qed
\end{prf}
 
The crucial result in \cite{BDLR,BG}  is the following  rigidity phenomenon of the angle between the minimal intermediate $C^*$-subalgebras.

\begin{thm}[\cite{BG}]\label{m1}(Rigidity of angle)
Let $\mathscr{B\subset A}$ be an irreducible inclusion of simple unital $C^*$-algebras with a conditional expectation $E:\mathscr{A\rightarrow B}$ of finite Watatani index. Then, the interior angle between any two distinct minimal intermediate $C^*$-subalgebras $\mathscr{C}$ and $\mathscr{D}$ of $\mathscr{B \subset A}$ is greater that $\pi/3$.
\end{thm}

\textbf{Proof of \Cref{bound}\,:} Following \cite{BG}, thanks to \Cref{m1}, we first note that $\texttt{\#}\,\mathcal{I}(\mathscr{B\subset A})\leq 9^{\text{dim}_{\mathbb{C}}(\mathscr{B}^{\prime}\cap \mathscr{A}_1)}$. In what follows, we prove that $\text{dim}_{\mathbb{C}}(\mathscr{B}^{\prime}\cap \mathscr{A}_1)\leq {[\mathscr{A:B}]}_0$.

Recall that given an inclusion of $C^*$-algebras $\mathscr{D\subset C}$ with a conditional expectation $E:\mathscr{C}\to\mathscr{D}$, the Pimsner-Popa probablistic constant \cite{PP}, denoted by $\lambda_E(\mathscr{D\subset C)}$, is defined as follows\,:
\[
\lambda_E(\mathscr{D\subset C})=\sup\{\lambda>0: E(x)\geq \lambda x~\text{for all}~x\in {\mathscr{C}}_{+}\}\,.
\]
Consider the following quadruple of $C^*$-algebras
\[
\begin{matrix}
\mathscr{A} &\subset & \mathscr{A}_1\\
\cup & &\cup \\
\mathscr{B}^\prime\cap\mathscr{A} =\bbc &\subset & \mathscr{B}^\prime\cap\mathscr{A}_1
\end{matrix}
\]
Denote by  $F$ the  conditional expectation from $\mathscr{B}^\prime\cap\mathscr{A}_1$ onto $\mathscr{B}^\prime\cap\mathscr{A}=\bbc$ obtained by  restricting the conditional expectation $E_1:\mathscr{A}_1\to\mathscr{A}$. This immediately says the following
\begin{eqnarray}\label{pp}
\lambda_{E_1}(\mathscr{A}\subset\mathscr{A}_1)\leq\lambda_F(\bbc\subset\mathscr{B}^\prime\cap\mathscr{A}_1).
\end{eqnarray}
In fact, since $E_0\circ E_1$ is a minimal conditional expectation, we see that $F$ is a trace (Thm. $2.12.3$, \cite{Wa}). Since $\mathscr{B}^\prime\cap\mathscr{A}_1$ is finite-dimensional $C^*$-algebra, assume that $\mathscr{B}^\prime\cap\mathscr{A}_1=\oplus_{j=1}^k\,M_{n_j}$, where $M_{n_j}$ denotes the algebra of $n_j\times n_j$ ($1\leq j\leq k$ for some $k\in\bbn$) matrices over $\bbc$. By $1.1.7$, Page $175$, \cite{Po}, (see also \cite{Hi} in this regard) we know that
\[
\lambda_F(\bbc\subset\mathscr{B}^\prime\cap\mathscr{A}_1)\leq\lambda_\tau(\bbc\subset\mathscr{B}^\prime\cap\mathscr{A}_1),
\]
where $\tau$ is the Markov trace for the inclusion $\bbc\subset\oplus_{j=1}^k\,M_{n_j}$. Therefore, \Cref{pp} gives us the following
\begin{eqnarray}\label{pq}
\lambda_{E_1}(\mathscr{A}\subset\mathscr{A}_1)\leq\lambda_\tau(\bbc\subset\oplus_{j=1}^k\,M_{n_j})\,.
\end{eqnarray}
We claim that the Markov trace $\tau$ on $\oplus_{j=1}^k\,M_{n_j}$ is the restriction of the unique normalized trace on the type $I$ factor $M_{n_1^2+\ldots+n_k^2}$, which happens to be the basic construction of $\bbc\subset\oplus_{j=1}^k\,M_{n_j}$. Suppose that $\textbf{t}=(t_1,\cdots, t_k)$ denotes the trace vector for $\oplus_{j=1}^k\,M_{n_j}$ corresponding to the Markov trace $\tau$. The inclusion matrix $\Lambda$ for $\bbc\subset\oplus_{j=1}^k\,M_{n_j}$ is simply the row matrix $(n_1\,\,n_2\,\cdots\,\,n_k)$ as the inclusion is unital, and we also have $\sum_{j=1}^kn_jt_j=1$. It is easy to verify that if we set $t_j=\frac{n_j}{n_1^2+\ldots+n_k^2}$ for each $j$, then the equation $\sum_{j=1}^kn_jt_j=1$ is satisfied, and moreover, $\Lambda^T\Lambda\textbf{t}=\alpha\textbf{t}$ is also satisfied for $\alpha=n_1^2+\ldots+n_k^2$. Thus, the required Markov trace $\tau$ on $\oplus_{j=1}^k\,M_{n_j}$ is given by the trace vector $\big(\frac{n_1}{n_1^2+\ldots+n_k^2},\ldots,\frac{n_k}{n_1^2+\ldots+n_k^2}\big)$. Observe that the constant $\alpha$ is the dimension of the type $I$ factor $M_{n_1^2+\ldots+n_k^2}$, and this matches with the Watatani index for the inclusion $\bbc\subset\oplus_{j=1}^k\,M_{n_j}$.

\noindent By Theorem $6.1$ in \cite{PP}, we have the following
\begin{eqnarray}\label{eqn}
\lambda_\tau(\bbc\subset\oplus_{j=1}^k\,M_{n_j})=\big\{\max_j\frac{1}{t_j}\big\}^{-1}=\big\{\max_j\frac{n_1^2+\ldots+n_k^2}{n_j}\big\}^{-1}=\frac{\min_j\,n_j}{n_1^2+\ldots+n_k^2}\,.
\end{eqnarray}
By \Cref{pq,eqn}, we get the following
\begin{IEEEeqnarray*}{lCl}
[\mathscr{A:B}]_0=\lambda_{E_1}(\mathscr{A}\subset\mathscr{A}_1)^{-1} &\geq& \lambda_\tau(\bbc\subset\oplus_{j=1}^k\,M_{n_j})^{-1}\\
&=& \frac{n_1^2+\ldots+n_k^2}{\min_j\,n_j}\\
&=& \frac{\mathrm{dim}\oplus_{j=1}^k\,M_{n_j}}{\min_j\,n_j}\\
&=& \frac{\mathrm{dim}(\mathscr{B}^\prime\cap\mathscr{A}_1)}{\min_j\,n_j}\,.
\end{IEEEeqnarray*}
Finally, since $\mathscr{B\subset A}$ is irreducible, using Theorem 4.6.3 in \cite{GHJ}, we get that $\min_j\,n_j=1$. Therefore, by \cite{BG} we finally have the following
\begin{center}
$\texttt{\#}\,\mathcal{I}(\mathscr{B\subset A})\leq 9^{\mathrm{dim}(\mathscr{B}^\prime\cap\mathscr{A}_1)}\leq 9^{[\mathscr{A:B}]_0}\,,$
\end{center}
which concludes the proof.\qed
\medskip

\noindent\textbf{Open Question:} Is it true that there are constants $s$ and $t$  such that for any irreducible inclusion of simple $C^*$ algebra $\mathscr{B\subset A}$ with a conditional expectation of index-finite type, the number of intermediate $C^*$-subalgebras is less than $s {[\mathcal{A:B}]_0}^t$?

\section{Acknowledgements}
We sincerely thank Professor Mikael Rørdam and Professor Yasuo Watatani for useful exchange. K. C. Bakshi acknowledges support of INSPIRE Faculty grant DST/INSPIRE/04/2019/002754 and S. Guin acknowledges support of SERB grant MTR/2021/000818.


\bigskip

\bigskip

\noindent {\em Department of Mathematics and Statistics,\\ Indian Institute of Technology Kanpur,\\ Uttar Pradesh $208016$, India}
\medskip

\noindent {Keshab Chandra Bakshi:} keshab@iitk.ac.in, bakshi209@gmail.com\\
{Satyajit Guin:} sguin@iitk.ac.in\\
{Debabrata Jana:} debabrata.jana05@gmail.com  


\bigskip

\begin{thebibliography}{10}

\bibitem{BDLR}
Bakshi, K. C.; Das, S.; Liu, Z.; Ren, Y.:
\newblock An angle between intermediate subfactors and its rigidity.
\newblock {\em Trans. Amer. Math. Soc.} {\bf 371} (2019), no.~8,  5973--5991.

\bibitem{BG}
Bakshi, K. C.; Gupta, V. P.:
\newblock Lattice of intermediate subalgebras.
\newblock {\em J. London Math. Soc.} {\bf 104(2)} (2021), 2082--2127.

\bibitem{BGS}
Bakshi, K. C.; Guin, S.; Sruthymurali:
\newblock Fourier-theoretic inequalities for inclusions of simple $C^*$-algebras.
\newblock {\em New York J. Math.} {\bf 29} (2023), 335--362.

\bibitem{GHJ}
Goodman, F.; de la Harpe, P.; Jones, V. F. R.:
\newblock Coxeter graphs and towers of algebras.
\newblock Mathematical Sciences Research Institute Publications, {\bf 14}. {\em Springer-Verlag, New York}, 1989.

\bibitem{Hi} Hiai, F.:
\newblock Minimizing indices of conditional expectations onto a subfactor. 
\newblock {\em Publ. Res. Inst. Math. Sci.} {\bf 24} (1988), 673--678.

\bibitem{Ino}
Ino, S.; Watatani, Y.:
\newblock Perturbation of intermediate $C^*$-subalgebras for simple $C^*$-algebras.
\newblock {\em Bull. Lond. Math. Soc.} {\bf 46} (2014), 469--480.

\bibitem{Iz}
Izumi, M.:
\newblock Inclusions of simple $C^*$-algebras.
\newblock {\em J. Reine Angew. Math.} {\bf 547} (2002), 97--138.

\bibitem{J}
Jones, V. F. R.:
\newblock Index for subfactors.
\newblock {\em Invent. Math.} {\bf 72} (1983), 1--25.

\bibitem{JS}
Jones, V. F. R.; Sunder, V. S.:
\newblock Introduction to subfactors.
\newblock London Mathematical Society Lecture Note Series, {\bf 234}. {\em Cambridge University Press, Cambridge}, 1997.

\bibitem{Ko}
Kosaki, H.:
\newblock Extension of Jones' theory on index to arbitrary factors.
\newblock {\em J. Funct. Anal.} {\bf 66} (1986), 123--140.

\bibitem{Lon}
Longo, R.:
\newblock Conformal subnets and intermediate subfactors.
\newblock {\em Comm. Math. Phys.} {\bf 237} (2003), 7--30.

\bibitem{PP}
Pimsner, M.; Popa, S.:
\newblock Entropy and index for subfactors.
\newblock {\em Ann. Sci. École Norm. Sup.} (4) {\bf 19} (1986), no. 1, 57--106.

\bibitem{Po}
Popa, S.:
\newblock Classification of amenable subfactors of type II.
\newblock {\em Acta Math.} {\bf 172} (1994), no. 2, 163--255.

\bibitem{PZ}
Pfender, F.; Ziegler, G.M.:
\newblock Kissing numbers, sphere packings and some unexpected proofs.
\newblock Notices Amer. Math. Soc. 51 (2004), 873--883.

\bibitem{R}
R{\o}rdam, M.:
\newblock Irreducible inclusions of simple $C^*$-algebras.
\newblock {\em Enseign. Math.} {\bf 69} (2023), no. 3/4, 275--314.

\bibitem{TW}
Teruya, T.; Watatani, Y.:
\newblock Lattices of intermediate subfactors for type $III$ factors.
\newblock {\em Arch. Math. (Basel)} {\bf 68} (1997), no. 6, 454--463.
  
\bibitem{Wa}
Watatani, Y.:
\newblock Index for $C^*$-subalgebras.
\newblock {\em Mem. Amer. Math. Soc.} {\bf 83} (1990), no. 424, vi+117 pp.

\bibitem{W}
Watatani, Y.:
\newblock Lattices of intermediate subfactors.
\newblock {\em J. Funct. Anal.} {\bf 140} (1996), no. 2, 312--334.

\bibitem{Z}
Zacharias, J.:
\newblock Splitting for subalgebras of tensor products.
\newblock {\em Proc. Amer. Math. Soc.} 129 (2001), no. 2, 407--413.

\bibitem{Zs}
Zsidó, L.:
\newblock A criterion for splitting $C^*$-algebras in tensor products.
\newblock {\em Proc. Amer. Math. Soc.} 128 (2000), no. 7, 2001--2006.

\end{thebibliography}
\end{document}